\documentclass{proc-l}
\usepackage{amscd}

\newtheorem{thm}{Theorem}[section]
\newtheorem{cor}[thm]{Corollary}
\newtheorem{lem}[thm]{Lemma}
\newtheorem{prop}[thm]{Proposition}

\newcommand{\thmref}[1]{Theorem~\ref{#1}}
\newcommand{\propref}[1]{Proposition~\ref{#1}}
\newcommand{\lemref}[1]{Lemma~\ref{#1}}

\def\hh{{\rm H}}
\def\dd{\partial}
\def\db{\bar \partial}

\def\noi{\noindent}
\def\ol{\overline}

\def\to{\rightarrow}
\def\too{\longrightarrow}
\def\zetaf{$\zeta$--function\ }

\def\CC{{\bf C}}
\def\RR{{\bf R}}
\def\ZZ{{\bf Z}}
\def\Aa{{\mathcal A}}
\def\Cc{{\mathcal C}}
\def\Dd{{\mathcal D}}
\def\Ee{{\mathcal E}}
\def\Ff{{\mathcal F}}
\def\Hh{{\mathcal H}}
\def\Kk{{\mathcal K}}
\def\Ll{{\mathcal L}}
\def\Ss{{\mathcal S}}
\def\Ww{{\mathcal W}}

\def\z{\zeta}

\def\Det{\mbox{\rm Det\,}}
\def\DET{\mbox{\rm DET\,}}
\def\dom{\mbox{\rm dom\,}}
\def\End{\mbox{\rm End}}
\def\Grass{\mbox{{$\mathcal Gr$}}}
\def\Hom{\mbox{\rm Hom}}

\def\Ind{\mbox{\rm Ind\,}}
\def\Ker{\mbox{\rm Ker}}

\def\ind{\mbox{\rm ind\,}}
\def\ran{\mbox{\rm range\,}}
\def\tr{\mbox{\rm tr\,}}
\def\Tr{\mbox{\rm Tr\,}}

\begin{document}

\title{Splitting the Curvature of the
 Determinant Line Bundle}

\author{Simon Scott}

\address{Department of Mathematics, King's College, Strand,
 London WC2R 2LS, United Kingdom.}

\email{sscott@mth.kcl.ac.uk}

\dedicatory{Dedicado a la memoria de Hugo Rojas 1973-1997}



\begin{abstract} It is shown that the determinant line bundle
associated to a family of Dirac operators over a closed
partitioned manifold $M=X^{0}\cup X^{1}$ has a canonical Hermitian
metric with compatible connection whose curvature satisfies an
additivity formula with contributions from the families of Dirac
operators over the two halves.
\end{abstract}

\maketitle

%% SECTION 1

\section{Introduction}

 Let $\Dd$ be a family of Dirac operators over a closed
even-dimensional manifold $M$ parameterized by a smooth manifold
$B$. Let $M$ be partitioned into codimension 0 submanifolds
$M=X^{0}\cup X^{1}$and let $\Dd^{0},\Dd^{1}$ be the families
obtained by restricting $\Dd$ to the two halves. The purpose of
this Letter is to communicate an additivity formula for the
curvature of the determinant line bundle $\Ll$ associated to $\Dd$
relative to the partition of $M$ and a choice of elliptic boundary
condition $P$ for the families $\Dd^{0},\Dd^{1}$. To do this, we
construct a Hermitian metric and compatible connection on $\Ll$
using a regularized trace and determinant.  Full details of the
proofs will appear in \cite{Sc97}.

The proof relies on a canonical splitting isomorphism defined for
each $P$ between $\Ll$ and the determinant line bundles associated
to the families $\Dd^{0}$ and $\Dd^{1}$. A similar but less
general splitting isomorphism has been proved by Piazza
\cite{Pi96} which leads in the limit to an additivity formula for
the $\zeta$-function curvature as the metric is blown-up in a
normal direction to $X^{0}\cap X^{1}$.  The additivity theorem we
prove here requires no asymptotic arguments, this is because the
connection we use is constructed directly from boundary data
arising from the partition, a different partition defines a
different connection. These results are concordant with principles
suggested by Topological QFT, and provide an explicit
representative for the local anomaly obeying the same splitting
principle as the Chern class of the determinant bundle.

There is a precise relation between the metric and connection on
$\Ll$ constructed here to the $\zeta$-function metric and
connection \cite{BiFr86}. The details of this are the subject of
an article currently in preparation. For preliminary results see
\cite{ScWo96} for the case of a self-adjoint elliptic boundary
value problem over an odd-dimensional manifold, and \cite{ScTo98,
BoScWo98} for a detailed proof of the relation between the
canonical and $\z$-function geometry in the case of a family of
elliptic boundary value problems in dimension one.

%% SECTION 2

\section{Grassmann Sections and the Determinant Bundle}

Let $\pi : Z\too B$ be a smooth fibration of manifolds with fibre
diffeomorphic to a closed connected even-dimensional manifold $M$.
The tangent bundle along the fibres $T(Z/B)$ is taken to be
oriented, spin and endowed with a Riemannian metric $g_{T(Z/B)}$.
Let $\Ss(Z/B)$ be the vertical spinor bundle and $E$ a Hermitian
coefficient bundle. Associated to this data one has a smooth
elliptic family of Dirac operators $\Dd=\{ D_{b}:b\in B
\}:\Hh\too\Hh$, where $\Hh = \pi_{*}(\Ss(Z/B)\otimes E)$ is the
infinite-dimensional Hermitian vector bundle on $B$ whose fibre at
$b$ is the Frechet space of smooth sections $\Hh_{b} =
C^{\infty}(M_{b},\Ss_{b}\otimes E_{b})$.
 The $Z_{2}$ bundle grading
$\Hh= \Ff^{+}\oplus\Ff^{-}$ into positive and negative chirality
fields defines families of chiral Dirac operators $\Dd^{\pm}:
\Ff^{\pm}\too\Ff^{\mp}$ and hence two families of
finite-dimensional vector spaces $\ker (\Dd^{\pm}) =\{\ker
(D_{b}^{\pm}):b\in B\}$. The Quillen determinant line bundle
$\Det\Ind\Dd^{+}$ is a complex line bundle over $B$ with fibre at
$b$ canonically isomorphic to the complex line $\wedge^{max}(\Ker
\Dd_{b}^{+})^{-1} \otimes \wedge^{max}\Ker (\Dd_{b}^{-})$ (see
\cite{BiFr86,Qu85}).

A connection on $\Hh_{Y}$ is defined as follows
\cite{Bi86,BiFr86}. A connection on the fibration means a
splitting
\begin{equation}\label{e:a}
TZ = T(Z/B) \oplus T^{H}Z
\end{equation}
and specifies an isomorphism $T^{H}Z\cong \pi^{*}TB$. A choice of
metric $g_{TB}$ on $TB$ determines a metric $g_{TZ} =
g_{T(Z/B)}\oplus \pi^{*}g_{TB}$ on $TZ$ with Levi-Civita
connection $\nabla^{TZ}$. We define a connection on $T(Z/B)$ by
$$\nabla^{T(Z/B)} = P_{Z/B}\nabla^{TZ}P_{Z/B} + c,$$ where
$P_{Z/B}$ is the orthogonal projection on $TZ$ with range $T(Z/B)$
and $c$ is the 1-form on $Z$ defined by
$[\frac{1}{2}i_{\xi}d(vol_{Z/B})]_{V} = c(\xi)vol_{Z/B}$. Here
$vol_{Z/B}$ is the volume form defined by $g_{T(Z/B)}$ regarded as
an n-form on $Z$ and $[.]_{V}$ means the vertical component. The
connection $\nabla^{T(Z/B)}$ lifts to a connection  on the
vertical spinor bundle, and endowing $E$ with a connection
compatible with its metric we obtain a connection
$\nabla^{\Ss(Z/B)}$ on each of the bundles $\Ss(Z/B)\otimes E$,
$\Ss^{\pm}(Z/B)\otimes E$. Because a section of the
infinite-dimensional bundle $\Hh$ is identified as a section of
$\Ss(Z/B)\otimes E$ we have a unitary connection $\nabla^{Z}$ on
$\Hh,\Ff^{\pm}$
\begin{equation}\label{e:b}
\nabla^{Z}_{\xi}s = \nabla^{\Ss(Z/B)}_{\xi_{H}}s,
\end{equation}
where $\xi_{H}$ is the horizontal lift of $\xi\in
C^{\infty}(B;TB)$. The correction term $c$  ensures that
$\nabla^{Z}$ is unitary.

By a splitting of the fibration $\pi : Z = Z^{0}\cup_{Y}Z^{1}\too
B$ we mean an embedded connected codimension 1 subfibration of
closed manifolds $\dd\pi : Y\too B$ defining on each fibre of $Z$
a partition $M_{b} = X^{0}_{b}\cup_{Y_{b}}X^{1}_{b}$. A splitting
defines families of chiral Dirac operators $\Dd^{+} =  \Dd^{0}\cup
\Dd^{1}$  with $\Dd^{i}= \{D^{i}_{b} = D^{+}_{b}|X_{b}^{i}:b\in
B\}$  associated to the fibrations  of manifolds with boundary
$\pi^{i}: Z^{i}\too B$ with $\dd Z^{i} = Y^{i}$. Here the
manifolds $Y^{i}$ are identified but such that $Y^{1}$ has the
reverse orientation to $Y^{0}=Y$. We assume that $g_{Z/B}$ splits
isometrically in a collar neighbourhood $V = [-1,1]\times Y$ of
$Y\equiv \{0\}\otimes Y$  as  $ du^{2} +  g_{Y/B}$, where $u$ is a
normal coordinate to $\dd Z$ and $g_{Y/B}$  the induced metric on
the boundary fibration, and that the Hermitian metrics on the
bundles $S(Z/B), E$ split similarly. In $V_{b}$ the operators
$D^{i}_{b}$ have the form $\sigma_{b}(\dd/\dd u +  A_{b})$ where
$\sigma_{b} : S_{Y_{b}} \to S_{Y_{b}}$ is a unitary bundle
isomorphism and $A_{b}$ is an elliptic self-adjoint operator
identified with the Dirac operator over $Y_{b}$. We make the
assumption that the operators of the boundary family $$\Dd_{Y}=\{
A_{b}:b\in B \}:\Hh_{Y}\too\Hh_{Y},$$ with $\Hh_{Y}$ the
infinite-dimensional bundle of boundary spinor fields, have
kernels of constant rank. (This is convenient, but not essential.)
From \eqref{e:a} there is an induced connection on the boundary
fibration
\begin{equation}\label{e:aa}
TY = T(Y/B) \oplus T^{H}Y
\end{equation}
and hence a connection $\nabla^{Y}$ on $\Hh_{Y}$ coinciding with
the analogue of the construction of $\nabla^{Z}$.

The bundle $\Hh_{Y}$ is endowed with a $\ZZ_{2}$-grading $\Hh_{Y}
= \Hh_{Y}^{+}\oplus\Hh_{Y}^{-}$ defined on each fibre by
$\Hh_{Y,b} = \hh_{b}^{+}\oplus \hh_{b}^{-}$ where $\hh_{b}^{+}$
(resp. $\hh_{b}^{-}$) is the direct sum of the eigenspaces of
$A_{b}$ with non-negative (resp. negative) eigenvalues. Associated
to the grading we have a (smooth) Grassmann bundle $\Grass
(\Hh_{Y})\too B$ with fibre the infinite-dimensional {\it smooth
Grassmannian}  $\Grass_{\infty}(\hh_{b})$ parameterising
projections $P\in\End (\hh_{b})$ which differ by a smoothing
operator from the projection $\Pi_{\geq,b}$ with image
$\hh_{b}^{+}$, where projection means self-adjoint indempotent.
(See \cite{BoWo93,PrSe90,Sc95,ScWo96}). A {\em Grassmann section}
for the family $\Dd^{0}$ is defined to be a smooth section $P$ of
the fibration $\Grass (\Hh_{Y})$. Reversing the orientation on $Y$
defines an opposite Grassmann bundle $ \Grass (\Hh_{Y^{1}})$, and
so, denoting the space of Grassmann sections by $\Grass (Y/B)$, we
have $P\in \Grass (Y/B)$ if and only if $I-P\in \Grass (Y^{1}/B)$.
In particular, since $\dim\ker A_{b}$ is assumed to be independent
of $b$, the spectral section $\Pi_{\geq}$ is an element of $\Grass
(Y/B)$. Thus any Grassmann section $P\in \Grass (Y/B)$ defines an
endomorphism of the bundle $\Hh_{Y}$ which differs from
$\Pi_{\geq}$ by a smooth family of smoothing operators.
\begin{prop}\label{p:a1}
A Grassmann section $P$ is equivalent to a smooth ungraded Frechet
subbundle $\Ww\too B$ of $\Hh_{Y}$. The bundle $\Ww$ is endowed
with a natural Hermitian metric with compatible connection
$\nabla^{\Ww} = P\cdot\nabla^{Y}\cdot P$ whose curvature is the
first-order differential operator
\begin{equation}\label{e:d}
\RR^{\Ww}(\xi_{1},\xi_{2}) = P\nabla^{\Ss
(Y/B)}_{-P_{Y/B}[\xi^{H}_{1},\xi^{H}_{2}]}\cdot P +
P\RR^{Y/B}(\xi^{H}_{1},\xi^{H}_{2})P
                  - S_{\Ww}^{*}\wedge S_{\Ww},
\end{equation}
where $\xi^{H}$ is the horizontal lift of $\xi\in
C^{\infty}(B;TB)$ and $\RR^{Y/B}$
 is the curvature of $\nabla^{T(Y/B)}$. Here $S_{\Ww} = \nabla^{Y}_{|\Ww} - \nabla^{\Ww}$ is
 the second fundamental form, given locally by
$S_{\Ww} = (I-P)\{dP +
\Theta\}P\in\Omega^{1}(U,\Hom(\Ww,\Ww^{\perp}),$ where $\Theta$ is
the connection 1-form for $\nabla^{Y}$ over $U$.
\end{prop}
\begin{proof}
 $\Ww$ is the bundle over $B$ with fibre $W_{b} =
\ran (P_{b})$. The first statement follows  from the local
triviality of  $\Hh_{Y}$ and hence of $\Grass (\Hh_{Y})$. The
induced metric is defined by $<s_{1},s_{2}>_{\Ww} =
<Ps_{1},Ps_{2}>$, where the right-side is the Hermitian metric on
$\Hh_{Y}$ and it is clear that the connection $\nabla^{\Ww}$ is
compatible. Equation \eqref{e:d} follows easily from \cite{Bi86}
Prop. 1.11. The local formula for $S_{\Ww}$ is an immediate
consequence of the definitions.
\end{proof}

The curvature of the connection $\nabla^{\Ww}$ is therefore not of
trace-class and so we cannot define a Chern character for $\Ww$.
There is, nevertheless, for each pair of Grassmann sections
$P^{0},P^{1}$ an associated determinant line bundle with
well-defined Chern class whose curvature is expressed as the
difference of the curvature operators of the associated bundles
$\Ww^{0},\Ww^{1}$. This is because (seen as operators on
$\Hh_{Y})$ they have the same top order symbol, and in fact differ
by only a trace-class operator.

To describe these constructions we use a Schatten class calculus,
as follows. Let $\End_{cpt}(\Ww)$ denote the algebra of compact
bundle endomorphisms (acting on the dense subbundle $\Ww$ of the
$L^{2}$ Hilbert bundle completion $\ol{\Ww}$). The subbundle
$\End_{tr}(\Ww)$ is the bundle whose fibre $\End_{tr}(W_{b})$ is
the ideal (first Schatten class) in $\End_{cpt}(W_{b})$ of
operators with  $\|A_{b}\|_{tr} =
\tr_{W_{b}}(A_{b}^{*}A_{b})^{1/2} < \infty,$ where the trace is
the sum over the eigenvalues.
 Each smooth section of $\End_{tr}(\Ww)$ has an associated trace function
 $b\mapsto \tr_{W_{b}}(A_{b})$.
More generally, there is a $\ZZ$-graded algebra of differential
forms taking values in $\End_{tr}(\Ww)$
$$\Omega^{k}(B,\End_{tr}(\Ww)) = \Omega^{k}(B)\otimes
\Omega^{0}(B,\End_{tr}(\Ww))$$ with multiplication
$\Omega^{i}_{tr} \times \Omega^{j}_{tr}\too\Omega^{i+j}_{tr}$ and
a smooth linear trace map $$\Tr_{\Ww}:
\Omega^{k}(B,\End_{tr}(\Ww))\too \Omega^{k}(B)\;\;\;
\theta_{b}\otimes A_{b}\mapsto \theta_{b}.\tr_{W_{b}}(A_{b}).$$
satisfying the supertrace identity
\begin{equation}\label{e:str}
\Tr_{\Ww}(A_{1}A_{2})  = (-1)^{({\rm deg}A_{1})({\rm
deg}A_{2})}\Tr_{\Ww}(A_{2}A_{1}).
\end{equation}
From the pointwise estimate $\|\wedge^{r} A_{b}\|_{tr}\leq
\frac{1}{r!}(\|A_{b}\|)^{r}$ we have  for each $r\in\ZZ_{+}$ a
well-defined algebra map
$$\Omega^{k}(B,\End_{tr}(\Ww))\too\Omega^{k}(B,\End_{tr}(\wedge^{r}\Ww))\;\;\;
\theta.A\mapsto\theta.\wedge^{r}A$$ and a smooth map, the Fredholm
determinant,
\begin{equation}\label{e:frdet}
\Omega^{0}(B,\End_{tr}(\Ww))\too \Omega^{0}(B),\;\;\;\; A_{b} \to
{\rm det}_{F}(1 + \epsilon A_{b}) :=
\sum_{r=0}^{\infty}\epsilon^{r}\tr_{\Ww_{b}}(\wedge^{r}A_{b}),
\end{equation}
The Fredholm determinant is multiplicative on elements in $1 +
\End_{tr}(\Ww)$, and characterized by the property:
\begin{lem}\label{l:2.1}
\begin{equation}\label{e:a2}
\frac{d}{d\epsilon}{\rm det}_{F}(1 + \epsilon A)_{|\epsilon=0}  =
\Tr_{\Ww}(A).
\end{equation}
\end{lem}

There are many different determinant line bundles to compare.
First, any $P\in \Grass(Y/B)$ defines a smooth family of elliptic
boundary value problems $$(\Dd^0,P)= \Dd^{0}: \dom (\Dd^{0},P)\too
 L^{2}(Z^{0},\Ss^{-}(Z/B)\otimes E)$$ $$\dom (\Dd^{0},P) =
\left\{\psi\in H^{1}(Z^{0},\Ss^{+}(Z/B)\otimes E):
P\psi_{|Y}=0\right\},$$ where $H^{1}$ is the first Sobolev space.
Similarly, we have $(\Dd^{1},I-P)$. The families $(\Dd^{0},P)$ and
$(\Dd^{1},I-P)$ define elliptic families of Fredholm operators
varying smoothly with the parameters, and we obtain
 the following fact.
\begin{prop}\label{p:2.5}
For each $P\in\Grass(Y/B)$ there is a Quillen  determinant line
bundle $\Det\Ind (\Dd^{0},P)$ with canonical determinant section
$b\mapsto \det D^{0}_{P_{b}}$ non-zero precisely where
$D^{0}_{P_{b}}$ is invertible. Similarly, we have $\Det\Ind
(\Dd^{1},I-P)$.
\end{prop}

Second, associated to each pair of Grassmann sections
$P^{0},P^{1}$ there is a canonical
 smooth Fredholm family of generalized Toeplitz operators
$$(P^{0},P^{1})\in C^{\infty}(B;\Hom (\Ww^{0},\Ww^{1})),\;\;\;
(P^{0},P^{1})_{b}\equiv P^{1}_{b}P^{0}_{b}: W_{0,b}\too W_{1,b},
$$ where $\Ww^{i}$ are the bundles of Proposition \ref{p:a1}. To
define the determinant bundle for such a family needs a more
general construction than the Quillen determinant line bundle.

To this end, we adapt a construction of Segal \cite{Se90}. By an
{\em admissible family} associated to the fibration $\pi : Z\to B$
we mean any of the families defined above or associated families
or any finite-rank homorphism of vector bundles over  $B$. Let
$\Aa = \{ A_{b}:b\in B\}:\Ee^{+}\too\Ee^{-}$ be such a family
where $\Ee^{\pm}$ are the appropriate bundles over $B$. The
determinant line $\Det (A_{b})$ when
 $\ind A_{b} = 0$ is defined to be the set of equivalences
classes $[T,\lambda ]$ of pairs $(T,\lambda )$ where $T- A_{b}\in
\Hom_{tr}(\Ee_{b}^{+},\Ee_{b}^{-})$ and $\lambda\in\CC$ relative
to the equivalence relation  $(Tq,\lambda ) \sim
(T,\det_{F}(q)\lambda )$ for $q\in 1 + \End_{tr}(\Ee_{b}^{+})$. If
$\ind A = r$  we define $\Det (A_{b})$ to be the determinant line
of  $A_{b}\oplus 0$ as an operator  $\Ee_{b}\too\Ee^{-} \oplus
\CC^{r}$ if  $r > 0$, or  $\Ee_{b} \oplus \CC^{-r}\too\Ee^{-}$ if
$r < 0$.

\begin{prop}\label{p:2.14}
Let $\Aa$ be an admissible family. Then $\DET (\Aa )= \cup_{b\in
B}\Det A_{b}$ is a complex line bundle over $B$ with canonical
section $\det$  non-zero if and only if  $A_{b}$ is invertible.
\end{prop}
\begin{proof}
To prove this we can work over the open covering of $B$ defined by
the sets
\begin{equation}\label{e:open}
U_{\alpha} = \{ b\in B : A_{b}+\alpha (b)\; {\rm invertible} \},
\end{equation}
where $\alpha$ is a finite-rank section  of $\Hom
(\Ee^{+},\Ee^{-})$. Complex multiplication is defined on $\Det
(A_{b})$ by $ \mu [T,\lambda ] = [T, \mu\lambda ], \;\;\mu\in\CC.$
For non-zero elements $[T_{\alpha (b)},\lambda ],
\newline [T_{\beta (b)},\lambda ]\in \Det A_{b}$ over $U_{\alpha}\cap U_{\beta}$, we have
$[T_{\alpha (b)},\lambda ] = g_{\alpha\beta}(A_{b})[T_{\beta
(b)},\lambda ]$ with $ g_{\alpha\beta}(A_{b}) = \det
_{F}((A_{b}+\alpha (b))(A_{b}+\beta (b))^{-1})$. Since
$(A_{b}+\alpha (b))(A_{b}+\beta (b))^{-1}=1+(\alpha (b)-\beta (b))
(A_{b}+\beta (b))^{-1})\in 1 + \End_{tr}(\Ww_{b})$, then from
\eqref{e:frdet} the map $b\mapsto g_{\alpha\beta}(A_{b})$ is a
well-defined smooth transition function on $U_{\alpha}\cap
U_{\beta}$. This defines $\DET (\Aa)$ globally as a complex line
bundle over B. The global canonical determinant section is $ \det:
B \too \DET (\Aa), \;\;\;b\too \det A_{b} \equiv [A_{b},1].$
\end{proof}

Thus, in addition to the determinant bundles $\DET (\Dd^{0},P)$,
we now have for each family $(P^{0},P^{1})$ a determinant line
bundle $\DET (P^{0},P^{1})$ over $B$. The following technical
lemma facilitates subsequent identifications.
\begin{lem}\label{lem:3.1}
Let $\Aa_{i} = \{A_{i,b}:\Ee^{+}_{i,b}\too\Ee^{-}_{i,b}\}$ for
$i=0,1,2$ be admissible families related by a commutative diagram
$$\begin{CD} 0 @>>> \Ee^{+}_{0}    @>>>  \Ee^{+}_{1}   @>>>
\Ee^{+}_{2}  @>>> 0 \\ @.   @VV{\Aa_{0}}V     @VV{\Aa_{1}}V
@VV{\Aa_{2}}V               \\ 0 @>>> \Ee^{-}_{0}     @>>>
\Ee^{-}_{1}   @>>>   \Ee^{-}_{2} @>>> 0.
\end{CD}$$
where the rows are short exact sequences of vector bundles,
 exact on each Sobolev
completion. Then there is a canonical isomorphism of determinant
line bundles
\begin{equation}\label{e:3.23}
\DET \Aa_{1} \cong \DET\Aa_{0} \otimes \DET\Aa_{2}
\end{equation}
which preserves the determinant sections: $\det (A_{1,b})
\longleftrightarrow \det (A_{0,b})\otimes \det (A_{2,b}).$
\end{lem}

In making the identifications of determinant bundles it is
convenient to work with the equivalent family of Dirac boundary
value problems
 $$\Dd^{0}\oplus P: H^{1}(Z^{0},\Ss^{+}(Z/B)\otimes E)\too L^{2}(Z^{0},\Ss^{-}(Z/B)\otimes E)
 \oplus \overline{\Ww}$$
$$(D^{0}\oplus P)(\psi) = (D^{0}\psi ,P(\psi_{|Y})).$$ The family
$\Dd^{0}\oplus P$ is elliptic with well-defined
 determinant line bundle $\Det\Ind(\Dd^{0}\oplus P).$
Using \lemref{lem:3.1} one can check that it is equivalent to work
with either of the constructions of the determinant line bundle
for the families $\Dd^{+},(\Dd^{0},P)$ and $(\Dd^{0}\oplus P)$ (
here $\Dd^{+}$ is the family of chiral Dirac operators over the
closed manifold $M$):

\begin{prop}\label{p:3.100}
There are canonical isomorphisms of determinant line bundles
preserving the determinant
 sections
\begin{equation}\label{e:3.300}
\Det\Ind \Dd^{+} \cong \DET \Dd^{+}
\end{equation}
\begin{equation}\label{e:3.31}
\Det\Ind (\Dd^{0},P) \cong \DET (\Dd^{0},P)\cong \Det\Ind
(\Dd^{0}\oplus P)\cong \DET(\Dd^{0}\oplus P).
\end{equation}
\end{prop}

The following relation is crucial in explaining the sewing
properties of determinant line bundles.

\begin{thm}\label{lem:3.2}
Let $P^{0},P^{1},P^{2}\in\Grass (Y/B) $. Then there is a canonical
isomorphism of determinant line bundles
\begin{equation}\label{e:3.30}
\DET (P^{0},P^{1})\otimes \DET (P^{1},P^{2})\cong \DET
(P^{0},P^{2})
 \end{equation}
with
\begin{equation}\label{e:3.301}
\det(P^{0},P^{1})\otimes  \det (P^{1},P^{2})\longmapsto
\det[(P^{1},P^{2})\circ (P^{0},P^{1})].
\end{equation}
\end{thm}

\begin{proof}
 Let $\Ww^{i}, i=0,1,2$ be the bundles defined by
the $P^{i}$ and let $U_{\alpha}$ be as in \eqref{e:open}. Over
$U_{\alpha}\cap U_{\beta}\neq\emptyset$ with
$\alpha\in\Hom(\Ww^{0},\Ww^{1}), \beta\in\Hom(\Ww^{1},\Ww^{2})$,
the isomorphism is defined by
\begin{equation}\label{e:3.103}
 \det(P^{0},P^{1})_{\alpha}\otimes
\det(P^{1},P^{2})_{\beta}\longmapsto \det
(P^{1},P^{2})_{\beta}\circ ((P^{0},P^{1})_{\alpha}.
\end{equation}
Because $(P^{1},P^{2})_{\beta}\circ (P^{0},P^{1})_{\alpha}$
differs from $(P^{0},P^{2})$ by the sum of a family of finite-rank
operators and a family of smoothing operators, and hence by an
family of operators of trace-class, we have $\DET (P^{0},P^{2}) =
\DET[(P^{1},P^{2})_{\beta}\circ (P^{0},P^{1})_{\alpha}].$ By an
obvious patching argument the proposition is proved.

Alternatively, the result is an immediate corollary of
\ref{lem:3.1} applied to the commutative diagram of bundle maps
$$
\begin{CD}
0 @>>> \Ww^{0}    @>{j_{1}}>>  \Ww^{0}\oplus \Ww^{1}
@>{j_{2}}>>    \Ww^{1}   @>>> 0 \\ @. @VV{(P^{0},P^{1})_{\alpha}}V
@VV{\Phi_{\alpha,\beta}}V @VV{(P^{1},P^{2})_{\beta}}V \\ 0 @>>>
\Ww^{1}     @>{j_{3}}>> \Ww^{2}\oplus \Ww^{1}  @>{j_{4}}>> \Ww^{2}
@>>> 0.
\end{CD}
$$
where $$\Phi_{\alpha,\beta} = (P^{1},P^{2})_{\beta}\circ
((P^{0},P^{1})_{\alpha}\oplus Id_{\Ww^{1}}$$ $$j_{1}(\xi ) =
(\xi,(P^{0},P^{1})_{\alpha}(\xi ))$$ $$j_{2}(\xi, \eta  ) =
(P^{0},P^{1})_{\alpha}(\xi ) - \eta$$ $$j_{3}(\xi ) =
((P^{1},P^{2})_{\beta}(\xi ),\xi)$$ $$j_{4}(\xi, \eta  ) = \xi
-(P^{1},P^{2})_{\beta}(\eta).$$
\end{proof}

%% SECTION 3

\section{Splitting Isomorphisms}

There is a canonical smooth Grassmann section which carries {\it
global} data about the family $\Dd^{0}$. This is the Calderon
section $P(\Dd^{0})$ defined fibrewise by the orthogonal
projection $P_{b}(\Dd^{0})$ onto the $L^{2}$ closure of Cauchy
data space of harmonic spinors $\Kk_{b}(\Dd^{0}) = \{ \phi\in
\Hh_{Y,b}: \phi = \psi_{|Y_{b}}, \psi\in \Ff^{+}_{b},
D^{0}_{b}\psi =0\}.$ That $P(\Dd^{0})$ is an element of $\Grass
(Y/B)$ depends on the metric $g_{Z/B}$ being a product in a
neighbourhood of the boundary $\dd Z$.

\begin{prop}\label{p:3.4}
Let $P^{0},P^{1},P$ be Grassmann sections for the family
$\Dd^{0}$. Then there is a canonical isomorphism of determinant
line bundles
\begin{equation}\label{e:3.1111}
\Det\Ind (\Dd^{0},P^{0}) \cong \Det\Ind(\Dd^{0}, P^{1})\otimes
\DET (P^{1},P^{0}),
\end{equation}
which preserves the canonical determinant sections if and only
$P^{0}=P^{1}$. One has
\begin{equation}\label{e:3.222}
\Det\Ind (\Dd^{0},P) \cong \DET (P(\Dd^{0}),P).
\end{equation}
\end{prop}

\begin{proof} Let $\Kk^{0},\Ww^{0}$ be the Frechet bundles over
$B$ defined by the Grassmann sections $P(\Dd^{0})$ and $P^{0}$.
Then we have the following commutative diagram of bundle maps with
exact rows $$\begin{CD} 0 @>>> \Kk^{0}    @>>>  \Ff^{+}_{0}
@>{\Dd^{0}}>>   \Ff^{-}_{0}  @>>> 0 \\ @.
@VV{(P^{0},P(\Dd^{0}))}V     @VV{\Dd^{0}\oplus P^{0}}V
@VV{Id}V               \\ 0 @>>> \Ww^{0}     @>>>
\Ff^{-}_{0}\oplus\Ww^{0}   @>>>   \Ff^{-}_{0} @>>> 0.
\end{CD}$$
The maps in the lower-row are the obvious ones. To see that the
upper-row is exact one needs to know that $\Dd^{0}$ defines a
surjective bundle map and that $\Kk^{0}$ is canonically identified
via the Poisson operator with the bundle with fibre $\Ker
D^{0}_{b}$ at $b\in B$. From Lemma  \ref{lem:3.1} there is a
canonical isomorphism of determinant line bundles $$\DET
(\Dd^{0}\oplus P^{0}) \cong \DET (P(\Dd^{0}), P^{0})\otimes \Det
Id\cong \DET (P(\Dd^{0}), P^{0}).$$ Combined with
\propref{p:3.100} this proves \eqref{e:3.222}, and
\eqref{e:3.1111} is now a consequence of \thmref{lem:3.2}.
\end{proof}

Notice that the analogue of \eqref{e:3.222} for the family
$\Dd^{1}$ is
\begin{equation}\label{e:opp}
\Det\Ind (\Dd^{1},I-P) \cong \DET (P,I-P(\Dd^{1})).
\end{equation}
A partition of the closed manifold $M=X^{0}\cup X^{1}$ induces a
splitting of the determinant line bundle in the following sense.

\begin{thm}\label{t:3.1}
For Grassmann sections $P^{0},P^{1}$ there is a canonical bundle
isomorphism
\begin{equation}\label{e:3.17}
\Det\Ind \Dd^{+} \cong \Det\Ind (\Dd^{0},P^{0})\otimes \Det\Ind
(\Dd^{1},I- P^{1}) \otimes \DET (P^{0},P^{1}).
\end{equation}
The determinant section $\det \Dd^{+}$ of $\Det\Ind \Dd^{+} $ maps
to the determinant section of the bundle on the right-side if and
only if $P^{0}=  P(\Dd^{0})$ and $P^{1}= I-P(\Dd^{1})$, in which
case
\begin{equation}\label{e:3.111}
\Det\Ind \Dd^{+} \cong \DET (P(\Dd^{0}), I-P(\Dd^{1})).
\end{equation}
\end{thm}

Note that the isomorphism \eqref{e:3.111} is completely intrinsic,
no choice of a Grassmann section is required. For the isomorphism
(\ref{e:3.17}) it is enough to prove (\ref{e:3.111}), the general
formula follows from Proposition \ref{p:3.4}. The proof of
(\ref{e:3.111}) uses Mayer-Vietoris type arguments, the details
will be presented in \cite{Sc97}, see \cite{Sc95} for the single
operator case and \cite{Se90}
 for the case of $\db$-operators over a Riemann surface.

As a corollary we obtain the following additivity formulae for the
Chern class of the determinant line bundle.

\begin{cor}\label{p:4.1}
 In $H^{2}(B)$
\begin{equation}\label{e:4.5}
c_{1}(\Det\Ind (\Dd^{0},P^{0})) = c_{1}(\Det\Ind (\Dd^{0},P^{1}))
+ c_{1}(\DET (P^{0},P^{1}))
\end{equation}
\begin{equation}\label{e:4.7}
c_{1}(\Det\Ind \Dd^{+} ) = c_{1}( \Det\Ind (\Dd^{0},P^{0})) +
c_{1}(\Det\Ind (\Dd^{1},I-P^{1})) + c_{1}( \DET (P^{0},P^{1}))
\end{equation}
\end{cor}

%% SECTION 4

\section{The canonical metric and connection}

We can do better however, and give a local version of the formulas
in Corollary \ref{p:4.1}. This requires the construction of a
Hermitian metric and compatible connection on the respective
determinant bundles using a regularized determinant and trace. One
way to do that is to use the \zetaf regularization (see Quillen,
Bismut and Freed \cite{Qu85,BiFr86} for closed manifolds, and
\cite{Pi96} for $b$-determinant bundles) which by recent results
of Wojciechowski \cite{Wo97} extends to $\Det\Ind(\Dd^{0},P)$.
Here, however, we use a canonical regularization scheme associated
to the splitting of the fibration and the Schatten class calculus.

For Grassmann sections $P,P^{0},P^{1}$ we have the families of
Laplacian operators $$\Delta^{+} =
(\Dd^{+})^{*}\Dd^{+},\;\;\;(\Delta^{0},P) =
 (\Dd^{0},P)^{*}(\Dd^{0},P),\;\;\;
\Delta_{(P^{0},P^{1})} = (P^{0},P^{1})^{*}(P^{0},P^{1}).$$
\begin{prop}\label{p:d.1}
Let $\Ww^{0}$ be the bundle over $B$ defined by $P^{0}$. Then
\begin{equation}\label{e:detcl}
\Delta_{(P^{0},P^{1})}- Id \in \Omega^{0}(B,\End_{tr}(\Ww^{0})).
\end{equation}
There are canonical identifications of determinant bundles
\begin{equation}\label{e:d2}
\Det\Ind (\Delta^{0},P) \cong \Det\Ind \Delta_{(P(D^{0}),P)}
\end{equation}
\begin{equation}\label{e:d3}
\Det\Ind \Delta^{+} \cong \Det\Ind \Delta_{(P(D^{0}),I-P(D^{1}))}
\end{equation}
with
\begin{equation}\label{e:det2det}
\det(\Delta^{0},P) \leftrightarrow
\det\Delta_{(P(D^{0}),P)}\;\;\;\; \det\Delta^{+} \leftrightarrow
\det\Delta_{(P(D^{0}),I-P(D^{1}))}
\end{equation}
\end{prop}
\begin{proof}
 Any two Grassmann sections differ by a smooth
family of smoothing operators and so the first assertion is
obvious.

We prove (\ref{e:d2}), the proof for (\ref{e:d3}) is similar.
\thmref{lem:3.2} says that associated to the Grassmann sections
$P(D^{0}),P$ there is a canonical isomorphism
$$\DET(P(D^{0}),P)\otimes \DET(P,P(D^{0})) \too
\DET\Delta_{(P(D^{0}),P)}$$ $$\det(P(D^{0}),P)\otimes
\det(P,P(D^{0})) \longmapsto \det[(P,P(D^{0}))(P(D^{0}),P)] =
\det\Delta_{(P(D^{0}),P)}.$$ From (\ref{e:3.222}) there is a
canonical isomorphism $\DET(P(D^{0}),P)\cong \Det\Ind (D^{0},P)$
preserving the determinant sections, and similarly, taking adjoint
diagrams, we obtain $\DET(P,P(D^{0}))\cong \Det\Ind (D^{0},P)^{*}$
with again $\det(P,P(D^{0}))\leftrightarrow \det(D^{0},P)^{*}$.
Since there is a canonical identification
$$\Det\Ind(D^{0},P)\otimes \Det\Ind(D^{0},P)^{*} \too
\Det\Ind(\Delta^{0},P)$$ which takes $\det(D^{0},P)\otimes
\det(D^{0},P)^{*}$ to $\det(\Delta^{0},P)$ the proposition is
proved.
\end{proof}

From \eqref{e:detcl} and \eqref{e:det2det} we can define, relative
to the splitting of $M$, a canonical regularization of the
Laplacian determinants by
\begin{equation}
{\rm det}_{\Cc}(\Delta^{0},P) = {\rm
det}_{F}\Delta_{(P(D^{0}),P)},\;\;\;\;\; {\rm det}_{\Cc}\Delta^{+}
= {\rm det}_{F}\Delta_{(P(D^{0}),I-P(D^{1}))}.
\end{equation}

\begin{thm}\label{t:4.2}
To each splitting $Z = Z^{0}\cup_{Y} Z^{1}$ of the fibration there
is a canonical Hermitian metric $\|\,.\,\|$ on $\Det\Ind \Dd^{+}$
(resp. on $\Det\Ind (\Dd^{0},P)$) with $$\| \det D^{+}_{b}\|^{2} =
{\rm det}_{{\mathcal C}}\Delta^{+}_{b},$$ if  $\Dd^{+}_{b}$ is
invertible and zero otherwise (resp. $$\| \det D^{0}_{P_{b}}\|^{2}
= {\rm det}_{{\mathcal C}}(\Delta^{0},P)_{b}$$ if
$\Dd^{0}_{P_{b}}$ is invertible and zero otherwise.)
\end{thm}
\begin{proof}
The complex line $\DET(P^{0},P^{1})_{b}$
 has an Hermitian inner-product defined  by
\begin{equation}\label{e:d5}
< [T_{b},\lambda ],  [\tilde{T}_{b},\mu]>= \ol{\lambda}\mu {\rm
det}_{F}(T_{b}^{*}\tilde{T}_{b}).
\end{equation}
 It is easy to check that this definition is independent of the choice of
 representative $(T_{b},\lambda)$ and hence that the inner-product is well defined.
Using the inner-product (\ref{e:d5}) a Hermitian  metric is
defined on $\DET(P^{0},P^{1})$ over  $U_{\alpha}$ by
$\|\det(P^{0},P^{1})_{\alpha}\|^{2}$.  It is immediate that this
is equal to $\det_{\Cc}\Delta_{(P^{0},P^{1})_{\alpha}}$. That the
locally defined metrics patch together to define a $C^{\infty}$
metric on $\DET(P^{0},P^{1})$ follows from the construction of the
determinant line and the
 transition functions $g_{\alpha\beta}$ of \propref{p:2.14}.
 The metric on $\Det\Ind (\Dd^{0},P)$ (resp. $\Det\Ind \Dd^{+}$) is the induced
 metric via the identifications of Proposition \ref{p:3.4}
(resp. Theorem \ref{t:3.1}).
\end{proof}

Next we construct a connection compatible with this metric using a
regularized trace. Let  $\nabla^{0},\nabla^{1}$ be the unitary
connections on the bundles $\Ww^{0},\Ww^{1}$ associated to
Grassmann sections $P^{0},P^{1}$, of \propref{p:a1}, and let
$\nabla^{Hom}$ be the induced connection on the bundle
$\Hom(\Ww^{0},\Ww^{1})$. The crucial fact about the connections
$\nabla^{0},\nabla^{1}$ is that though they do not have local
trace-class connection forms, the connection they induce on
$\End(\Ww^{0})$ does.
\begin{prop}\label{p:d6}
Over $U_{\alpha}$ one has
\begin{equation}\label{e:d.11}
    (P^{0},P^{1})_{\alpha}^{-1}\nabla^{Hom}(P^{0},P^{1})_{\alpha}\in
\Omega^{1}(U_{\alpha},\End_{tr}(\Ww^{0})).
\end{equation}
The locally defined connections
\begin{equation}
\nabla^{\alpha}\det(P^{0},P^{1})_{\alpha} =
\Tr_{\Ww^{0}}\left[(P^{0},P^{1})_{\alpha}^{-1}\nabla^{Hom}(P^{0},P^{1})_{\alpha}\right]
\det(P^{0},P^{1})_{\alpha}
\end{equation}
patch together into a connection $\nabla^{Det}$ on
$\DET(P^{0},P^{1})$ which is unitary for the metric $\|.\|$.
\end{prop}
\begin{proof}
 We may assume that $\Hh_{Y}$ has been trivialized
over the open set $U_{\alpha}$, so that $\nabla^{Y} = d + \Theta$
for some 1-form $\Theta\in\Omega^{1}(U_{\alpha},\End(\Hh_{Y}))$.
We compute
$$(P^{0},P^{1})_{\alpha}^{-1}\nabla^{Hom}(P^{0},P^{1})_{\alpha}$$
$$   =  (P^{0},P^{1})_{\alpha}^{-1}\left\{d(P^{0},P^{1})_{\alpha}
+ P^{1}\Theta P^{1}(P^{0},P^{1})_{\alpha} -
(P^{0},P^{1})_{\alpha}P^{0}\Theta P^{0}\right\}.$$

Since $P^{1} = P^{0} + \Ss$ for some smooth family of smoothing
operators
 $\Ss = \{S_{b}\}:\Hh_{Y}\too\Hh_{Y}$,
the bracketed term is equal to $$ d(P^{0},P^{1})_{\alpha} +
P^{1}\Theta P^{1}P^{0} + P^{1}\Theta P^{1}\alpha - P^{1}P^{0}
\Theta P^{1} -  \alpha P^{0}\Theta P^{0} $$ $$ =  \; dP^{0} +
dP^{0}\Ss  + P^{0}\Ss  + d\alpha  + P^{0}\Theta \Ss P^{0} +
\Ss\Theta P^{0} - \Ss\Theta \Ss $$ $$\;\;\; - \Ss P^{0}\Theta
P^{0} +  P^{0}\Theta P^{0}\alpha + P^{0}\Theta \Ss + \Ss\Theta
P^{0}\alpha + \Ss\Theta \Ss\alpha,$$ which is a sum of smoothing
and finite-rank bundle operators, while
$$((P^{0},P^{1})_{\alpha})^{-1} = (P^{0} + P^{0}\Ss +
\alpha)^{-1}$$ and $\Ss$ is smoothing and $\alpha$ finite-rank.
Hence (\ref{e:d.11}) is a sum of smoothing and finite-rank
families and hence of trace-class.

The second and third statements are consequences of the identities
over \newline $U_{\alpha}\cap U_{\beta}$ $$ d\,{\rm
det}_{F}((P^{0},P^{1})_{\beta}^{-1}(P^{0},P^{1})_{\alpha}   = $$
$$ \left\{
\Tr_{\Ww^{0}}((P^{0},P^{1})_{\alpha}^{-1}\nabla^{Hom}(P^{0},P^{1})_{\alpha})
 -   \Tr_{\Ww^{0}}((P^{0},P^{1})_{\beta}^{-1}\nabla^{Hom}(P^{0},P^{1})_{\beta} \right\}$$
\begin{equation}\label{e:id1}
\times\; {\rm
det}_{F}((P^{0},P^{1})_{\beta}^{-1}(P^{0},P^{1})_{\alpha}
\end{equation}
and $$ d\,{\rm
det}_{F}((P^{1},P^{0})_{\beta}^{*}(P^{0},P^{1})_{\alpha}  = $$
$$\left\{\Tr_{\Ww^{0}}((P^{0},P^{1})_{\alpha}^{-1}\nabla^{Hom}(P^{0},P^{1})_{\alpha})
 -  \Tr_{\Ww^{0}}((P^{1},P^{0})_{\beta}^{-1}\nabla^{Hom}(P^{1},P^{0})_{\beta}\right\}$$
\begin{equation}\label{e:id2}
 \times \; {\rm det}_{F}((P^{1},P^{0})_{\beta}^{*}(P^{0},P^{1})_{\alpha},
\end{equation}
which utilize Lemma \ref{l:2.1}.
\end{proof}

Let $\nabla$ be the induced unitary connection on $\Hom
(\Ff^{+},\Ff^{-})$ from $\nabla^{Y}$. We define a connection  on
$\Det\Ind\Dd^{+}$ over $U_{\alpha}$  by
\begin{equation}\label{e:d12}
\nabla^{(\Dd^{+})}_{\alpha}\det \Dd^{+}_{\alpha} = \Tr_{\mathcal
C} \left[(\Dd_{\alpha}^{+})^{-1}\nabla\Dd_{\alpha}^{+}\right]\det
\Dd^{+}_{\alpha}.
\end{equation}
Here $\det(\Dd^{+}_{\alpha})$ is the image of $\det
(P(\Dd^{0}),I-P(\Dd^{1}))_{\alpha}$ under the canonical
isomorphism \eqref{e:3.111} (so $\det\Dd^{+}_{\alpha=0} =
\det\Dd^{+}$) and the connection 1-form is the regularized {\em
canonical trace} defined by $$\Tr_{\mathcal
C}\left[(\Dd_{\alpha}^{+})^{-1}\nabla\Dd_{\alpha}^{+}\right] :=
\Tr_{\Kk^{0}}\left[(P(\Dd^{0}),I-P(\Dd^{1}))_{\alpha}^{-1}\nabla^{Hom}
(P(\Dd^{0}),I-P(\Dd^{1})_{\alpha}\right]$$ where $\Kk^{0}$ is the
bundle over $B$ defined by the Calderon section $P(\Dd^{0})$.

Similarly, a connection can be defined on $\Det\Ind(\Dd^{0},P)$
over $U_{\alpha}$  by
\begin{equation}\label{e:d13}
\nabla_{\alpha}^{(\Dd^{0},P)}\det (D^{0},P)_{\alpha} =
\Tr_{\mathcal
C}\left[(\Dd^{0},P)_{\alpha}^{-1}\nabla(\Dd^{0},P)_{\alpha}\right]
\det (D^{0},P)_{\alpha}
\end{equation}
where $\det(\Dd^{0},P)_{\alpha})$ is the image of $\det
(P(\Dd^{0}),P)_{\alpha}$ under the canonical isomorphism
\eqref{e:3.222} (so $\det (\Dd^{0},P)_{\alpha=0} =
\det(\Dd^{0},P)$) and the connection 1-form is the regularized
(canonical) trace defined by
$$\Tr_{\mathcal
C}\left[(\Dd^{0},P)_{\alpha}^{-1}\nabla(\Dd^{0},P)_{\alpha}\right]
:= \Tr_{\Kk^{0}}\left[(P(\Dd^{0}),P)_{\alpha}^{-1}\nabla^{Hom}
(P(\Dd^{0}),P)_{\alpha}\right].$$ The 1-forms act as endomorphisms
of the determinant line bundles $\Det\Ind\Dd^{+}_{|U_{\alpha}}$
and $\Det\Ind(\Dd^{0},P)_{|U_{\alpha}}$ via the identifications of
Section 3.

\noi We then have:

\begin{thm}\label{t:d2}
For each splitting $Z = Z^{0}\cup_{Y} Z^{1}$ of the fibration the
connections
 defined locally over $U_{\alpha}$ by (\ref{e:d12})
(resp. (\ref{e:d13}) ) patch together into a connection
$\nabla^{(\Dd^{+})}$ on
 $\Det\Ind \Dd^{+}$ (resp. $\nabla^{(\Dd^{0},P)}$ on $\Det\Ind (\Dd^{0},P)$)
 compatible with the Hermitian metric.
\end{thm}

%% SECTION 5

\section{Additivity of the curvature}

So far, we know that associated to the fibration of closed
manifolds $Z\too B$ there is the  determinant line bundle
$\Det\Ind \Dd^{+}$ and that when endowed with a splitting
$Z^{0}\cup_{Y}Z^{1}$ there is a connection $\nabla^{(\Dd^{+})}$
canonically
 defined on $\Det\Ind \Dd^{+}$. If we make a choice of Grassmann section
 $P\in \Grass (Y/B)$ then we have determinant line bundles
$\Det\Ind (\Dd^{0},P)$ and $\Det\Ind (\Dd^{1},I- P)$ with unitary
connections $\nabla^{(\Dd^{0},P)}$ and $\nabla^{(\Dd^{1},I-P)}$.
\begin{thm}\label{t:add}
Let $\RR^{(\Dd^{+})}$ be the curvature 2-form of the connection
$\nabla^{(\Dd^{+})}$, and let
$\RR^{(\Dd^{0},P)},\RR^{(\Dd^{1},I-P)}$ be the curvatures of the
connections $\nabla^{(\Dd^{0},P)}$ and $\nabla^{(\Dd^{1},I-P)})$.
Then
\begin{equation}\label{e:split}
\RR^{(\Dd^{+})} = \RR^{(\Dd^{0},P)} + \RR^{(\Dd^{1},I-P)}.
\end{equation}
\end{thm}
Note that \eqref{e:split} is taking place in $\Omega^{2}(B)$,
since the endomorphism bundle of a complex line bundle is
canonically trivial.
\begin{proof} From Theorem \ref{t:3.1} there is an isomorphism
\begin{equation}\label{e:isom}
\Det\Ind \Dd^{+} \stackrel{\pi}{\too} \Det\Ind (\Dd^{0},P)\otimes
\Det\Ind (\Dd^{1},I- P)
\end{equation}
and hence connections $\nabla^{(\Dd^{+})}$ and
$\nabla^{\pi}=\pi^{*}(\nabla^{(\Dd^{0},P)}\otimes 1 + 1\otimes
\nabla^{(\Dd^{1},I-P)})$ canonically defined on the determinant
line bundle $\Det\Ind\Dd^{+}$. The proof of the theorem proceeds
over the open sets where $(\Dd^{0},P)_{\alpha}$ and
$(\Dd^{1},I-P)_{\beta}$ are invertible. For notational brevity we
just do the case over the region $U$ where $\alpha,\beta =0$, note
that   $\Dd^{+}$ is invertible over $U$. The relation between the
corresponding connection 1-forms $\omega^{+} = \Tr_{{\mathcal
C}}(\Dd^{+})^{-1}\nabla\Dd^{+},$ $\;\omega^{0}_{P} =
\Tr_{{\mathcal C}}(\Dd^{0},P)^{-1}\nabla(\Dd^{0},P),$ and \newline
$\omega^{1}_{I-P} =  \Tr_{{\mathcal
C}}(\Dd^{1},I-P)^{-1}\nabla(\Dd^{1},I-P)$ can be computed by
evaluating $\nabla^{(\Dd^{+})}$ and $\nabla^{\pi}$ on the section
 $\det\Dd^{+}$. Using the identifications of \thmref{lem:3.2} and Sect.3 we find that
$$\omega^{+} = \omega^{0}_{P} + \omega^{1}_{I-P} + d\,\log F,$$
where $F$ is the function on $U$ defined by $$F(b) = {\rm
det}_{F}\left[\frac{(P(\Dd_{b}^{0}),I-P(\Dd_{b}^{1})}
{(P(\Dd_{b}^{0}),P_{b})(P_{b},I-P(\Dd_{b}^{1})}\right].$$
Since the connection forms differ by only a closed 1-form, the
theorem is proved.
\end{proof}

The additivity may also be deduced from the following identities
expressing the curvature of the determinant line bundles in terms
of the curvatures of the bundles defined by the corresponding
Grassmann sections.
\begin{thm}\label{t:families}
Over $U_{\alpha}$ one has:
\begin{equation}\label{e:R1}
\RR^{(\Dd^{+})} =
\Tr_{\Kk^{0}}\left[(P(\Dd^{0}),I-P(\Dd^{1}))^{-1}\RR^{\Kk^{1}}(P(\Dd^{0}),I-P(\Dd^{1}))
- \RR^{\Kk^{0}}\right]
\end{equation}
\begin{equation}\label{e:R2}
\RR^{(\Dd^{0},P)} =
\Tr_{\Kk^{0}}\left[(P(\Dd^{0}),P)^{-1}\RR^{\Ww}(P(\Dd^{0}),P) -
\RR^{\Kk^{0}}\right]
\end{equation}
\begin{equation}\label{e:R3}
\RR^{(\Dd^{1},I-P)} =
\Tr_{\Ww}\left[(P,I-P(\Dd^{1}))^{-1}\RR^{\Kk^{1}}(P,I-P(\Dd^{1}))
- \RR^{\Ww}\right]
\end{equation}
Here $\RR^{\Kk^{0}},\RR^{\Kk^{1}},\RR^{\Ww}$ are the curvature
forms of \propref{p:a1} for the bundles $\Kk^{0}, \Kk^{1}, \Ww$
defined by the Grassmann sections $P(\Dd^{0}),I-P(\Dd^{1}),P$.
\end{thm}
\thmref{t:families} may be construed as a statement of the local
families index theorem, expressing the curvature of the
determinant line bundle as the degree 2 component of a local
families index density.
\begin{proof} The essential point is that
although the curvature form $\RR^{\Ww}$ defined by a Grassmann
section $P$ is not of trace-class, the differences on the
right-hand side of \eqref{e:R1},\eqref{e:R2},\eqref{e:R3} are of
trace-class. The formulas follow by a direct computation using the
local connection forms in \eqref{e:d12} and \eqref{e:d13}.
\end{proof}

The additivity  \eqref{e:split} now follows from the identities
$$\Tr_{\Ww^{0}}\left[\RR^{\Ww^{0}} -
(P^{0},P^{1})^{-1}\RR^{\Ww^{1}}(P^{0},P^{1})\right] =
\Tr_{\Ww^{1}}\left[(P^{0},P^{1})\RR^{\Ww^{0}}(P^{0},P^{1})^{-1} -
\RR^{\Ww^{1}}\right]$$
$$\Tr_{\Ww^{0}}\left[(P^{0},P^{2})^{-1}\RR^{\Ww^{2}}(P^{0},P^{2})
- \RR^{\Ww^{0}}\right] = \hspace{6.3cm}$$ $$\hspace{4.0cm}
\Tr_{\Ww^{0}}\left[[(P^{1},P^{2})(P^{0},P^{1})]^{-1}\RR^{\Ww^{2}}[(P^{1},P^{2})(P^{0},P^{1})]
- \RR^{\Ww^{0}}\right]$$ where $P^{i}\in \Grass (Y/B)$ define
bundles $\Ww^{i}$ with curvature $\RR^{\Ww{i}}$, $i=0,1,2$. The
first follows by writing the operators with respect to the direct
sum bundle $\Ww^{0}\oplus\Ww^{1}$ and using \eqref{e:str}. The
second, since the terms differ by an identically vanishing trace
over an anti-commutator term.

The additivity is particularly simple in the case where the
boundary fibration $\dd Z=Y\too B$ is trivial and $\nabla^{Y}$ is
the trivial connection, since everything in sight is then
trace-class. The connection on $\Ww$ defined by a Grassmann
section $P$ is then $\nabla^{\Ww} = P\cdot d\cdot P.$ From
\eqref{e:d} or by direct computation, this has trace-class
curvature form with $\Tr_{\Ww}\left[\RR^{\Ww}\right] =
\Tr_{\Ww}(PdPdP)$. \thmref{t:families} then reduces to:
\begin{thm}\label{t:trivY}
\begin{equation}
\RR^{(\Dd^{+})} = \Tr_{\Kk^{1}}\left[\RR^{\Kk^{1}}\right] -
\Tr_{\Kk^{0}}\left[\RR^{\Kk^{0}}\right]
\end{equation}
\begin{equation}
\RR^{(\Dd^{0},P)} =  \Tr_{\Ww}\left[\RR^{\Ww}\right] -
\Tr_{\Kk^{0}}\left[\RR^{\Kk^{0}}\right]
\end{equation}
\begin{equation}
\RR^{(\Dd^{1},I-P)} = \Tr_{\Kk^{1}}\left[\RR^{\Kk^{1}}\right] -
\Tr_{\Ww}\left[\RR^{\Ww}\right].
\end{equation}
\end{thm}

%% BIBLIOGRAPHY

\end{document}